\documentclass[12pt]{article}
\begin{document}
\title{A problem of Stanis\l{}aw Saks}
\author{Alexandre Eremenko\thanks{Supported by NSF grant DMS-1665115.}}
\date{}
\maketitle
\begin{abstract}
A solution of Problem 184 from the Scottish Book is given.

2010 MSC 31A05. Keywords: subharmonic functions.
\end{abstract}

On February 8, 1940, the following entry was made
in the Scottish book~\cite{S}:
\vspace{.1in}

{\em
\noindent
184. Problem; S. Saks.
\newline
A subharmonic function $\phi$ has everywhere partial derivatives
$\partial^2\phi/\partial x^2,\;\partial^2\phi/\partial y^2$.
Is it true that $\Delta\phi\geq 0$?
\newline
Remark: it is obvious immediately that $\Delta\phi\geq 0$ at all
points of continuity of 
$\partial^2\phi/\partial x^2,\;\partial^2\phi/\partial y^2$,
therefore on an everywhere dense set.
\newline
Prize: one kilo of bacon}.
\vspace{.1in}

\noindent
{\bf Theorem.}
{\em Let $u$ be a subharmonic function of two variables 
whose first partial derivatives exist on the coordinate axes
and $u_{xx},\; u_{yy}$ exist at the origin.
Then 
$u_{xx}(0,0)+u_{yy}(0,0)\geq 0.$
}
\vspace{.1in}

{\em Proof.} 
Without loss of generality we 
assume that $u(0,0)=u_x(0,0)=u_y(0,0)=0$ (add a linear function).
Proving the Theorem by contradiction, we assume that
$\Delta u(0,0)<0$. Then there exist real $a,b$ and $R_0>0$ such that
for $x^2+y^2<R_0^2$ we have
\begin{equation}\label{1}
u(x,0)\leq a x^2,\quad u(0,y)\leq by^2,\quad
\mbox{where}\quad a+b<0.
\end{equation}
Without loss of generality, $a<0$.

If $b< 0$, consider the function
$$v_1(r\cos\theta,r\sin\theta)=Cr^2|\sin(2\theta)|,$$
which is harmonic in each quadrant,
and choose $C>0$ so large that
$v_1(x,y)\geq u(x,y)$ when $x^2+y^2=R_0^2$.
Then $u(x,y)\leq v_1(x,y)$ for $x^2+y^2<R_0^2$ by the
Maximum principle applied to the intersection of this disk 
with each quadrant.
Thus 
\begin{equation}\label{3}
u(x,y)\leq C(x^2+y^2),\quad\mbox{when}\quad x^2+y^2<R_0^2.
\end{equation}

Consider the family of subharmonic functions
$$u_r(x,y)=r^{-2}u(rx,ry),\quad r>0$$
In view of (\ref{3}), for every compact $K$ in the plane there exists
$r_0>0$ such that $u_r$ are defined and uniformly bounded from above
on $K$ for $r\in(0,r_0)$.
Therefore there is a sequence $r_j\to 0$
for which 
$u_{r_j}\to u_0$ in $L^1_{\mathrm{loc}}$, 
where $u_0$ is a subharmonic function, \cite[Theorem 3.2.12]{H}.
Moreover
$u(x,y)\geq\limsup_{r\to 0} u_0(x,y)$ for every $x,y$
by \cite[Theorem 3.2.13]{H}, so $u_0(0,0)=0$.
To show that $u_0$ satisfies (\ref{1}), fix a point $(x_0,0)$,
and consider disks $B_t$ of radii $t$ centered at this point. Since
the family $\{ u_r\}$ is uniformly bounded from above on $B_1$,
there is a continuous majorant $v$ for this family in $B_1$,
such that $v(x_0,0)\leq ax_0^2$. This $v$ is just the solution
of the Dirichlet problem for upper and lower halves of $B_1$
with boundary conditions $ax^2$ on the intersection of $B_1$
with the $x$-axis, and
constant on the half-circles. So for every $\epsilon>0$ there exists
$\delta$ such that $v(x_0,0)\leq ax_0^2+\epsilon$ in $B_\delta$.
Then $L^1_{\mathrm{loc}}$ convergence gives
$$u_0(x_0,0)\leq
\frac{1}{|B_\delta|}\int_{B_\delta}u_0(x,y)dxdy\leq
\frac{1}{|B_\delta|}\int_{B_\delta}v(x,y)dxdy\leq
ax_0+\epsilon.$$
As $\epsilon$ is arbitrary, we obtain that $u_0$
satisfies the first inequality in (\ref{1}) on the whole $x$-axis.
Similar arguments show that $u_0$ satisfies the
second inequality in (\ref{1}) on the whole $y$-axis,
and also satisfies (\ref{3}) in the whole plane.
 
The Phragm\'en--Lindel\"of indicator of $u_0$,
$$h(\theta):=\limsup_{r\to\infty}r^{-2}u_0(r\cos\theta,r\sin\theta)$$
is non-positive for $\theta=\pi/2$
and negative for $\theta=0$.
This contradicts the inequality
$$h(\theta)+h(\theta+\pi/2)\geq 0,$$
which the indicators of all functions of order $2$ must satisfy,
\cite[Section 8.2.4]{L}.

If $b\geq 0$, we consider the subharmonic
function
$$u^*(x,y)=u(x,y)+c(x^2-y^2),$$
where $b<c<-a$. Such a $c$ exists because $a+b<0$ in (\ref{1}). 
Then $u^*$ satisfies
$$u^*(x,0)\leq (a+c)x^2,\quad u^*(0,y)\leq (b-c)y^2$$
near the origin, and
we apply the previous argument to $u^*$.
This completes the proof.
\vspace{.1in}

\noindent
{\em Corollary.} There is no subharmonic function $u$ satisfying
$$u(0)=0\quad\mbox{and}\quad u(x,0)\leq-\epsilon|x|$$
for all sufficiently small $x$ and  $\epsilon>0$. 
\vspace{.1in}

\noindent
{\em Remark.} The Theorem does not hold in $R^n$ for $n\geq 3$.
Indeed, in this case the union of the coordinate axes is a polar set,
so it is easy to construct a counterexample.

{\em Purdue University,

West Lafayette, IN 47907

USA

eremenko@purdue.edu}
\end{document}